\documentclass[12pt]{amsart}
\usepackage{amsmath,amsfonts,euscript,amscd,amsthm,amssymb,upref,graphics}

\theoremstyle{plain}
\newtheorem{theorem}{Theorem}
\swapnumbers

\newtheorem{proposition}[subsection]{Proposition}
\newtheorem{lemma}[subsection]{Lemma}
\newtheorem{corollary}[subsection]{Corollary}

\newtheorem{convention}[subsection]{Convention}

\theoremstyle{definition}
\newtheorem{definition}[subsection]{Definition}
\newtheorem{example}[subsection]{Example}
\newtheorem{remark}[subsection]{Remark}

\newtheorem{nothing*}[subsection]{}
\newtheorem{notation}[subsection]{Notation}

\newcommand{\rien}[1]{}

\newcommand{\AVF}{ \operatorname{{\rm AVF}}}
\newcommand{\VFH}{ \operatorname{{\rm VF}_{hol}}}

\newcommand{\LieH}{ \operatorname{{\rm Lie}_{hol}}}
\newcommand{\VFA}{ \operatorname{{\rm VF}_{alg}}}
\newcommand{\LieA}{ \operatorname{{\rm Lie}_{alg}}}

\newcommand{\Aut}{ \operatorname{{\rm Aut}}}

\newcommand{\C}{\ensuremath{\mathbb{C}}}

\newcommand{\tX}{{\tilde X}}

\newcommand{\tZ}{{\tilde Z}}

\newcommand{\tS}{{\tilde S}}

\newcommand{\tq}{{\tilde q}}
\newcommand{\tr}{{\tilde r}}
\newcommand{\trho}{{\tilde \rho}}

\newcommand{\sgoth}{{\ensuremath{\mathfrak{s}}}}
\newcommand{\lgoth}{{\ensuremath{\mathfrak{l}}}}

\newcommand{\mgoth}{{\ensuremath{\mathfrak{m}}}}

\newcommand{\ggoth}{{\ensuremath{\mathfrak{g}}}}

\newcommand{\cL}{{\ensuremath{\mathcal{L}}}}
\newcommand{\cT}{{\ensuremath{\mathcal{T}}}}

\newcommand{\cE}{{\ensuremath{\mathcal{E}}}}

\newcommand{\cC}{{\ensuremath{\mathcal{C}}}}
\newcommand{\cH}{{\ensuremath{\mathcal{H}}}}
\newcommand{\Ker}{{\rm Ker} \,}

\renewcommand{\epsilon}{\varepsilon}
\renewcommand{\phi}{\varphi}
\renewcommand{\emptyset}{\varnothing}
\begin{document}
\renewcommand{\baselinestretch}{1.07}

\title[Criteria for the density property of complex manifolds ] {Criteria for the density property of complex manifolds}
\author{Shulim Kaliman}
\address{Department of Mathematics\\
University of Miami\\
Coral Gables, FL 33124 \ \ USA}
\email{kaliman@math.miami.edu}
\author{Frank Kutzschebauch}
\address{Mathematisches Institut \\Universit\"at Bern
    \\Sidlerstr. 5
   \\ CH-3012 Bern, Switzerland}
\email{Frank.Kutzschebauch@math.unibe.ch}
\thanks{{\bf Acknowledgments:} This research was started during a visit
of the second author to the University of Miami, Coral Gables, and
continued during a visit of both of us to the Max Planck Institute
of Mathematics in Bonn. We thank these institutions for their
generous support and excellent working conditions. The research of
the first author was also partially supported by NSA Grant no.
H98230-06-1-0063 and the second one by Schweizerische
Nationalfonds grant No 200021-107477/1. }
{\renewcommand{\thefootnote}{} }
%
\maketitle \vfuzz=2pt

\vfuzz=2pt
\section{Introduction}
The ground-breaking papers of Anders\' en and Lempert (\cite{A},
\cite{AL}) established remarkable properties of the automorphism group
of $\C^n$ $(n \ge 2 )$ which imply, in particular, that any local
holomorphic phase flow on a Runge domain $\Omega$ in $\C^n$ can be
approximated by global holomorphic automorphisms of $\C^n$ (for an
exact statement see Theorem 2.1 in \cite{FR}).

The next step in the development of the Anders\'en-Lempert theory
was made by Varolin who extended it from Euclidean spaces to a
wider class of algebraic complex manifolds. He realized also that
the following density property is crucial for this theory.
\begin{definition}\label{dense}

A complex manifold $X$ has the density property if in the
compact-open topology the Lie algebra $\LieH (X)$ generated by
completely integrable holomorphic vector fields on $X$ is dense in
the Lie algebra $\VFH (X)$ of all holomorphic vector fields on
$X$. An affine algebraic manifold $X$ has the algebraic density
property if the Lie algebra $\LieA (X)$ generated by completely
integrable algebraic vector fields on it coincides with  the Lie
algebra $\VFA (X)$ of all algebraic vector fields on it (clearly,
the algebraic density property implies the density property).

\end{definition}

In this terminology the main observation of  the Anders\'en-Lempert theory says that
$\C^n$ ($n\ge 2$) has the algebraic density property. Varolin and
Toth (\cite{V1}, \cite {TV1}, \cite{TV2}) established the density
property for some manifolds including semi-simple complex Lie
groups and some homogenous spaces of semi-simple Lie groups. Their
proof relies heavily on representation theory and does not, for example, lead to an answer in
the case of  other linear
algebraic groups.

In this paper we suggest new effective criteria for the density
property. This enables us to give a trivial proof of the original
Anders\'en-Lempert result and to establish (almost free of charge)
the algebraic density property for all linear algebraic groups
whose connected components are different from tori or $\C_+$. As
another application of this approach we tackle the question (asked
among others by F. Forstneri\v c) about the density of algebraic
vector fields on Euclidean space vanishing on a codimension 2
subvariety.

Our method of establishing the algebraic density property for an
affine algebraic variety $X$ consists of two ingredients described
in section 2. First, we try to find a nontrivial $\C [X]$-module
$L$ (over the algebra $\C [X]$ of regular functions on $X$) in
$\LieA (X)$. It turns out that this requires the existence of two
commuting completely integrable algebraic vector fields on $X$
satisfying some compatibility condition (see, Definition
\ref{compatible} below). Second, since $\LieA (X)$ is invariant
under algebraic automorphisms of $X$, in the presence of some
homogeneity property of $X$ we can increase $L$ so that it
coincides with the $\C [X]$-module of all algebraic vector fields
(in which $L$ is contained, of course, as a submodule). In
sections 3 and 4 we develop technique for checking this
compatibility condition and apply it in the cases of linear
algebraic groups and the complements to codimension 2 subvarieties
in Euclidean spaces.

{\em Acknowledgments.} We would like to thank D. Akhiezer for
inspiring discussions and consultations, F. Donzelli for catching
some inaccuracies, and the referee for valuable comments that
lead, in particular, to the present formulation of Proposition
\ref{condition2}.


\section{New approach to the Anders\'en-Lempert theory.}

The homogeneity property mentioned before is reflected in the
following.

\begin{definition}\label{tangential} Let $X$ be
an algebraic manifold and  $x_0\in X$. A finite subset $M$ of the
tangent space $T_{x_0} X$ is called a generating set if the image
of $M$ under the action of the isotropy subgroup of $x_0$ (in the
group of all algebraic automorphisms $\Aut X$ of $X$) generates
the whole space $T_{x_0}X$.

The manifold $X$ will be called tangentially semi-homogeneous if
it is homogeneous (with respect to $\Aut X$) and admits a
generating set consisting of one vector.
\end{definition}

\begin{theorem}\label{semi} Let $X$ be a homogeneous (with respect to $\Aut X$) affine
algebraic manifold with algebra of regular functions
$\C [X]$, and  $L$ be a submodule  of the $\C [X]$-module of all
vector fields such that $L \subset \LieA (X)$. Suppose that the
fiber of $L$ over some $x_0\in X$ contains a generating set. Then
$X$ has the algebraic density property.

\end{theorem}

\begin{proof} The $\C [X]$-modules $TX$ and $L$ generate coherent
sheaves $\cT$ and $\cL$ on $X$ where $\cL$ is a subsheaf of $\cT$.
The action of $\alpha \in \Aut X$ maps $\cL$ onto another coherent
subsheaf $\cL_{\alpha}$ of $\cT$. The sum of such subsheaves with
$\alpha$ running over a finite subset of $\Aut X$ is a coherent
subsheaf $\cE$ of $TX$. Let $\mgoth$ be the maximal ideal for
$x_0$. Definition \ref{tangential} implies that $\cE$ can be
chosen so that $\cE / \mgoth \cE$ coincides with $T_{x_0}X$.
Furthermore, since $X$ is homogeneous we can suppose that this is
true for every point in $X$. Thus $\cE = \cT$ (\cite{Har}, Chapter
II, exercise 5.8). Since composition with automorphisms preserves complete integrability, all global sections of $\cE$ are in
$\LieA (X)$ which concludes the proof.

\end{proof}

Another ingredient of our method is rooted in a new proof of the
following fact.

\begin{corollary}\label{AL} {\rm (The main observation of the Anders\'en-Lempert theory)}
For $n \geq 2$ the space $\C^n$ has the algebraic density property.

\end{corollary}

\begin{proof} Let $x_1, \ldots ,x_n$ be a coordinate system on
$\C^n$ and $\delta_i =\partial /\partial x_i$ be the partial
derivative, i.e. $\Ker \delta_i$ is the ring of polynomials
independent of $x_i$. Hence the polynomial ring $\C^{[n]}$ is
generated as a vector space by elements of $\Ker \delta_1 \cdot
\Ker \delta_2$. Note also that for $f_i \in \Ker \delta_i$ the
algebraic vector fields $f_i \delta_i$ and $x_if_i\delta_i$ are
completely integrable. Then the field
$$[f_1\delta_1,x_1f_2\delta_2]-[x_1f_1\delta_1,f_2\delta_2]=f_1f_2\delta_2$$
belongs to $\LieA (X)$ since $x_1f_2 \in \Ker \delta_2$. Thus
$\LieA (X)$ contains all algebraic fields proportional to
$\delta_2$. Since $\C^n$ is clearly tangentially semi-homogeneous
Theorem \ref{semi} implies the desired conclusion.
\end{proof}

\begin{remark} There is no need to  use tangential
semi-homogeneity in this proof since we can replace $\delta_2$ by
any other partial derivative $\delta_i$ and obtain each algebraic
vector field as a sum of fields proportional to $\delta_i, \, i=1,
\ldots , n$. However there are no partial derivatives on affine
algebraic varieties different from Euclidean spaces. Hence
formalization of this argument requires Theorem \ref{semi} and
some substitution for partial derivatives.

We remind that for a
completely integrable algebraic vector field on an affine
algebraic variety its phase flow is only a holomorphic
$\C_+$-action that is not necessarily algebraic.
\end{remark}

\begin{definition} An algebraic vector field $\delta$ on $X$ is called
semi-simple if its phase flow generates an algebraic $\C^*$-action
on $X$. A vector field $\sigma$ is called locally nilpotent if its
phase flow is an algebraic $\C_+$-action on $X$. In the last case
$\sigma$ can be viewed as a locally nilpotent derivation on the
algebra $\C [X ]$ of regular functions on $X$. That is, for every
$f \in \C [X]$ there is $n=n(f)$ for which $\sigma^n (f)=0$.

\end{definition}

There is one-to-one correspondence between the set of locally
nilpotent derivations on $\C [X]$ and the set of algebraic
$\C_+$-actions on $X$  (e.g., see \cite{Fre} for details). We
shall often use below the fact that for any locally nilpotent
derivation $\sigma$ and a regular function $f$ from its kernel
$\Ker \sigma$ (resp. regular function $f$ of degree 1 with respect
to $\sigma$, i.e. $\sigma (f) \in \Ker \sigma \setminus 0$) the
vector field $f \sigma$ is locally nilpotent (resp. completely
integrable).

\begin{definition}\label{compatible}
Let $\delta_1$ and $\delta_2$ be nontrivial algebraic vector
fields on an affine algebraic manifold $X$ such that $\delta_1$ is
a locally nilpotent derivation on $\C [X]$, and $\delta_2$ is
either also locally nilpotent or semi-simple. That is, $\delta_i$
generates an algebraic action of $H_i$ on $X$ where $H_1 \simeq
\C_+$ and $H_2$ is either $\C_+$ or $\C^*$. We say that $\delta_1$
and $\delta_2$ are compatible if (i) the vector space ${\rm Span}
(\Ker \delta_1 \cdot \Ker \delta_2)$ generated by elements from
$\Ker \delta_1 \cdot \Ker \delta_2$ contains a nonzero ideal in
$\C [X]$ and (ii) some element $a \in \Ker \delta_2$ is of degree
1 with respect to $\delta_1$, i.e. $\delta_1 (a) \in \Ker \delta_1
\setminus \{ 0 \}$.

\end{definition}

\begin{remark}\label{instead} Instead of condition (ii) suppose now
that $\delta_1$ and $\delta_2$ commute. Then by Corollary
\ref{non-degenerate} below condition (i) implies that the
$H_1$-action on $X$ generates a nontrivial algebraic $\C_+$-action
on $X//H_2$. Taking an element of $\C [X// H_2]$ whose degree with
respect to this action is 1 we can treat its lift-up to $X$ as
$a$. That is, for such commutative $\delta_1$ and $\delta_2$
condition (ii) is automatic.
\end{remark}

\begin{theorem}\label{density}
Let $X$ be a smooth homogeneous affine algebraic manifold with
finitely many pairs of compatible vector fields $\{ \delta_1^k
, \delta_2^k \}_{k=1}^m$ such that for some point $x_0 \in  X$
vectors $\{ \delta_2^k (x_0) \}_{k=1}^m$ form a generating set.
Then $\LieA (X)$ contains a nontrivial $\C [X]$-module and $X$ has
the algebraic density property.

\end{theorem}

\begin{proof}
Let $\delta_1$ and $\delta_2$ be one of our pairs.
Choose an element $a\in \Ker \delta_2$ of degree 1 with respect to
$\delta_1$ and set $b= \delta_1 (a)$. Let $f_i \in \Ker \delta_i$.
Then $[af_1\delta_1, f_2\delta_2 ] - [f_1\delta_1, af_2\delta_2]
=-b f_1f_2 \delta_2$. The last vector field is from $\LieA (X)$
and since $\delta_1$ and $\delta_2$ are compatible, Definition
\ref{compatible} implies that sums of such vector fields include
every vector field of form $I \delta_2$ where $I$ is a nonzero
ideal in $\C [X]$. Applying this argument to all compatible pairs
we see that $\LieA (X)$ contains all linear combinations of
$\delta_2^k$ with coefficients in some nonzero ideal $J\subset \C
[X]$. Since under a small perturbation of $x_0$ the set  $\{
\delta_2^k (x_0) \}_{k=1}^m$ remains a generating set we can
suppose that $x_0$ does not belong to the zero locus of $J$. Hence
by Theorem \ref{semi} $X$ has the algebraic density property.
\end{proof}

\begin{remark}
If $X$ is tangentially semi-homogenous and, furthermore, any non-zero tangent vector (at any point) is
a generating set, then  Theorem \ref{density}  implies that for the algebraic density property a single pair of
compatible vector fields is enough.
\end{remark}

\begin{corollary}\label{product} Let $X_1$ and $X_2$ be homogeneous affine
algebraic  manifolds such that each  $X_i$ admits a
finite number of integrable algebraic vector fields $\{ \delta_i^k
\}_{k=1}^{m_i}$ whose values at some point $x_i\in X_i$ form a
generating set and, furthermore, in the case of $X_1$ these vector
fields are locally  nilpotent. Then $X_1 \times X_2$ has the
algebraic density property.
\end{corollary}

\begin{proof} Note that $\delta_1^k$ and $\delta_2^j$ generate
compatible integrable vector fields on $X_1\times X_2$ which we
denote by the same symbols. Applying isotropy groups one can
suppose that $\{ \delta_i^k(x_i) \}$ is a basis of $T_{x_i}X_i$.
In order to show that the set of  vectors  $M=\{0\times
\delta_2^k(x_2) \}$ form a generating set in $T_{x_1\times x_2}
(X_1 \times X_2)$ we need the following fact that is obvious in a
local coordinate system.

{\bf Claim.} {\em Let $X$ be a complex manifold and let $\nu$ be a
vector field on $X$. Suppose that $f$ is a holomorphic function
from $\Ker \nu$ and $x_0 \in f^{-1}(0)$. Then phase flow induced
by the vector field $f\nu$ generates a linear action on the
tangent space $T_{x_0}X$ given by the formula $w \to w + df (w )
\nu (x_0)$ where $df$ is the differential and $w \in T_{x_0} X$.
In particular, the span of the orbit of $w$ under this phase flow
contains vector $df (w ) \nu (x_0)$.}

Applying this claim for $\nu = \delta_1^j$ we see that the orbit of $M$ under the isotropy group
of $x_1 \times x_2$ contains all vectors of form $\delta_1^j (x_1) \times \delta_2^k(x_2)$. Thus
$M$ is a generating set and we are done by Theorem \ref{density}.

\end{proof}

\begin{remark}
The reason why we use the locally nilpotent $\delta_1^j$ in this proof as $\nu$ and not
(the possibly semi-simple) $\delta_2^j$ is the following: The vector field $f\delta_2^j$ with $f \in \Ker\delta_2^j$
may not generate an algebraic action while $f\delta_1^j$ with $f \in \Ker\delta_1^j$ always generates an algebraic action.
It is worth mentioning if one wants to prove density property instead of algebraic density property the use
of $\delta_2^j$ is permissible.
\end{remark}


\begin{example}\label{ex3}

(1) Let $X=\C^k\times (\C^*)^l$ with $k\geq 1$ and $k+l\geq 2$.
Then $X$ has algebraic density property by Corollary \ref{product}.

(2) If $G$ is a simple Lie group then it is tangentially
semi-homogeneous since the adjoint action of $G$ generates an
irreducible representation on the tangent space $\ggoth$ at the
identity $e$ (i.e., any nonzero vector in $T_eG$ is a generating
set). Let $X$ be $SL_n (\C )$ with $n\geq 2$, i.e. $X$ is
tangentially semi-homogeneous. Then every $x \in X$ is a matrix
$(c_{kj})$ with determinant 1. Set $\delta_1 (c_{1j}) =c_{nj}$ and
$\delta_1 (c_{kj})=0$ for $k \ne 1$.  Set $\delta_2 (c_{nj})
=c_{1j}$ and $\delta_2 (c_{kj})=0$ for $k \ne n$. 
Note that constants and functions depending on $c_{kj}, k \ne 1$
are in $\Ker \delta_1$ while constants and functions depending on
$c_{kj}, k \ne n$ are in $\delta$. Therefore, condition (i) of
Definition \ref{compatible} holds. Taking $c_{11}$ as $a$ in
condition (ii) we see that $\delta_1$ and $\delta_2$ are
compatible. Thus $SL_n (\C )$ has the algebraic density property.

\end{example}

\section{Density of Affine Algebraic Groups with Connected Components
Different from Tori or $\C_+$.}

We start with a digest of the notion of categorical (algebraic)
quotient and its properties which will be used extensively in the
rest of this section.

\begin{definition}
Let $G$ be a linear algebraic group acting algebraically on an
affine algebraic variety $X$ and, therefore, on its algebra $\C
[X]$ of regular functions (we are going to use this notation for
the algebra of regular functions further for any variety $X$, not
necessarily affine). Consider the subalgebra $\C [X]^G$ of
$G$-invariant functions. Its spectrum is called the categorical
quotient of this action and it is denoted by $X//G$. The
monomorphism $\C [X]^G \hookrightarrow \C [X]$ generates a
dominant (but not necessarily surjective) morphism $\rho : X \to
X//G$ which is called the quotient morphism. The universal
property of categorical quotients says that any morphism from $X$
that is constant on orbits of $G$ factors through $\rho$.

\end{definition}

\begin{remark}\label{dimension}
For a reductive $G$ the subalgebra $\C [X]^G$ is always finitely
generated by Nagata's theorem and, therefore, $X//G$ can be viewed
as an affine algebraic variety. Furthermore, $\rho$
is surjective in this case, the points of $X//G$ are in one-to-one
correspondence with closed orbits of $G$ in $X$, and every fiber
of $\rho$ is the union of those orbits whose closure contains the
corresponding closed orbit (e.g., see \cite{Sch}). In particular,
if each orbit is closed then the categorical quotient coincides
with the usual geometric quotient (this happens, say, when a
reductive subgroup acts on a linear algebraic group by
multiplication). If $G$ is not reductive then $\C [X]^G$ is not
finitely generated in general (by Nagata's counterexample to the
fourteenth Hilbert problem). However, $X//G$ can be viewed as a
quasi-affine algebraic variety and $\C [X]^G$ as its algebra of
regular functions \cite{Wi03}. We will work mostly with $G \simeq
\C_+$. In this case general fibers of $\rho$ are always orbits of
the $\C_+$-action (i.e. $\dim X//G = \dim X -1$) and $\C[X]^G$
coincides with the kernel of the corresponding locally nilpotent
derivation (e.g., see \cite{Fre}).
\end{remark}

\begin{notation}\label{H_i} In this section
$H_1$ is isomorphic to $\C_+$ and $H_2$ is isomorphic either to
$\C_+$ or $\C^*$. We suppose also that $X$ is a normal affine
algebraic variety equipped with nontrivial algebraic $H_i$-actions
where $i=1,2$ (in particular, each $H_i$ generates an algebraic
vector field $\delta_i$ on $X$). The categorical quotients will be
denoted $X_i=X//H_i$ and the quotient morphisms by $\rho_i : X \to
X_i$.
\end{notation}

We start with a geometric reformulation of Definition
\ref{compatible}.

\begin{proposition}\label{reformulation} Set $\rho = (\rho_1, \rho_2) : X
\to Y:= X_1 \times X_2$ and $Z$ equal to the closure of $\rho (X)$
in $Y$.  Then ${\rm Span} (\Ker \delta_1 \cdot \Ker \delta_2)$
contains an ideal of $\C [X]$ iff $\rho : X \to Z$ is a finite
birational morphism.

\end{proposition}

\begin{proof} Every nonzero element of ${\rm Span} (\Ker \delta_1
\cdot \Ker \delta_2)$ is of the form $g \circ \rho$ where $g\in \C
[Z]=\C [Y]|_Z$. Thus ${\rm Span} (\Ker \delta_1 \cdot \Ker
\delta_2)$ coincides with the subalgebra $\rho^* (\C [Z ]) \subset
\C [X ]$ and we need to establish when $\rho^* (\C [Z ])$ contains
a nontrivial ideal of $\C [X]$. Note that functions from any
nontrivial ideal separate general points of $X$ while functions
from $\rho^* (\C [Z ])$ do not separate points of $\rho^{-1}(z)$
for any $z \in Z$. Hence $\rho : X \to Z$ must be birational if we
want $\delta_1$ and $\delta_2$ to be compatible.

Assume now that the closure of $Z \setminus \rho (X)$ contains a
divisor $D \subset Z$. There is a rational function $f$ on $Z$ so
that it has poles on $D$ and nowhere else. Multiplying $f$ by $h
\in \C [Z]$ such that $h$ is not identically zero on $D$ but
vanishes on $D \cap \rho (X)$ with sufficient multiplicity, one
can suppose that $f \circ \rho$ is regular on $X$. On the other
hand for $n$ sufficiently large and $g$ as before $gf^n$ has poles
on $D$ and cannot be a regular function on $Z$. Thus $(gf^n) \circ
\rho \notin {\rm Span} (\Ker \delta_1 \cdot \Ker \delta_2)$ and
the last vector space cannot contain a nonzero ideal in this case.

Let $\rho =\nu \circ \rho_0$ where $\nu : Z_0 \to Z$ is a
normalization of $Z$ and $\rho_0 : X \to Z_0$ is the induced
morphism of normal varieties. Then $\rho_0$ is birational and $Z_0
\setminus \rho_0 (X)$ is of codimension at least 2 since otherwise
even $\rho_0^*(\C [Z_0])$ does not contain a nontrivial ideal of
$\C [X]$. The indeterminacy set $V\subset \rho_0 (X)$ of the
rational map $\rho_0^{-1}$ is of codimension at least 2. Hence any
regular function on $\rho_0 (X)\setminus V$ extends to $Z_0$ by
the Hartogs theorem. This implies that $\rho_0^{-1}$ is regular
and, therefore, $\rho_0$ is an isomorphism, i.e. $\rho$ is finite
birational. Since $\nu$ is finite $\C [Z_0]$ is generated over $\C
[Z]$ by a finite number of functions of form $f_i/g_i, \, i=1,2,
\ldots ,n$ where $f_i$ and $g_i$ are regular on $Z$. Treat $\C [Z]
\simeq \nu^* (\C [Z])$ as a subalgebra of $\C [Z_0]$ and consider
the the principal ideal $J$ in $\C [Z_0]$ generated by
$\prod_{i=1}^n g_i$. By construction, $J\subset \C [Z]$. Hence
$\rho^* (\C [Z])$ contains a nonzero ideal of $ \C [X]\simeq \C
[Z_0]$ which is the desired conclusion.

\end{proof}

Note that for every (resp. a general) $(x_1,x_2) \in Z$ the set
$\rho^{-1}(x_1,x_2)=\rho_1^{-1}(x_1) \cap \rho_2^{-1}(x_2)$ is
finite (resp. a singleton) in this Proposition. Hence a
non-constant orbit of $H_1$ cannot be contained in a fiber of
$\rho_2$ and we have the following.

\begin{corollary}\label{non-degenerate}
In the case of $[\delta_1, \delta_2]=0$ the $H_1$-action on $X$
generates a nontrivial $\C_+$-action on $X_2$.

\end{corollary}

\begin{lemma}\label{finereformulation}
Suppose that $X,H_i,X_i, \delta_i$, and $\rho_i$ are as in
Notation \ref{H_i}, and either {\rm (i)} $[\delta_1 ,
\delta_2]=0$; or {\rm (ii)} $\delta_1$ and $\delta_2$ are both
locally nilpotent and generate a Lie algebra $\sgoth \lgoth_2$
that induces an algebraic action of $SL_2 (\C)$ on $X$. Set
$\Gamma = H_1 \times H_2$ in case (i), and $\Gamma = SL_2 (\C )$
in case (ii). Suppose that $X'$ is a normal affine algebraic
variety equipped with a non-degenerate (meaning that the generic orbits have dimension $2$ resp. $3$, equal to the dimension of $\Gamma$) $\Gamma$-action and $p : X
\to X'$ is a finite $\Gamma$-equivariant morphism (for each
$i=1,2$), i.e. we have commutative diagrams
\[ \begin{array}{ccc}
X & \stackrel{\rho_i}{\rightarrow} & X_i\\ \, \, \, \downarrow p&&
\, \,  \, \, \downarrow q_i \\ X' &
\stackrel{\rho_i'}{\rightarrow} &
X_i'\\
\end{array} \]
\noindent where $\rho_i' : X' \to X_i'=X'//H_i$ is the quotient
morphism of the $H_i$-action on  $X'$ (i. e., we treat $\C [X_i']$
as a subalgebra of $\C [X']$). Let ${\rm Span} (\C [X_1] \cdot \C
[X_2])$ contain a nonzero ideal of $\C[X]$. Then ${\rm Span} (\C
[X_1'] \cdot \C [X_2'])$ contains a nonzero ideal of $\C [X']$.

\end{lemma}

\begin{proof}
Since $p$ is finite, every $f \in \C [X_i] \subset \C [X]$ is a
root of a minimal monic polynomial with coefficients in $\C [X' ]$
that are constant on $H_i$-orbits (since otherwise $f$ is not
constant on these orbits). By the universal property these
coefficients are regular on $X_i'$, i.e. $f$ is integral over $\C
[X_i']$ and $q_i$ is finite. Consider the commutative diagram
\[ \begin{array}{ccc}
X & \stackrel{\rho}{\rightarrow} & X_1\times X_2\\ \, \, \, \,
\downarrow p && \, \,  \, \, \downarrow q \\ X' &
\stackrel{\rho'}{\rightarrow} &
X_1'\times X_2'\\
\end{array} \]
\noindent where $\rho =(\rho_1,\rho_2), \,  \rho'
=(\rho_1',\rho_2')$, and $q=(q_1,q_2)$. Let $Z$ (resp. $Z'$) be
the closure of $\rho (X)$ in $X_1\times X_2$ (resp. $\rho' (X')$
in $X_1'\times X_2'$). By Proposition \ref{reformulation} $\rho
(X)=Z$ and, therefore, (since $q$ is finite) $q (\rho (X))=\rho'
(X')=Z'$. Let $\nu : Z_0 \to Z$ be a normalization, i.e. $X$ is
naturally isomorphic to $Z_0$. Since  $q \circ \nu : Z_0 \to Z'$
is finite it generates a finite morphism $X \simeq Z_0 \to Z_0'$
onto a normalization $Z_0'$ of $Z'$. The commutativity of our
diagram implies that $\rho'$ generates a finite morphism $\rho_0'
: X' \to Z_0'$. Thus it suffices to prove the following.

{\em Claim.} In the last commutative diagram of
$\Gamma$-equivariant morphisms the fact that $\rho : X \to Z$ is
birational and $p$ (and, therefore, $q$) is finite implies that
morphism $\rho' : X' \to Z'$ is birational.

For any $x\in X$ we set $x'=p(x)$, $x_j=\rho_j (x)$, and
$x_{j}'=q_j(x_j)=\rho_j' (x')$. Assume that $x$ is a general point
of $X$ and $y \in X$ is such that $\rho' (x')=
(x_{1}',x_{2}')=(y_{1}',y_{2}')=\rho' (y')$. Hence $q_j^{-
1}(x_{j}')=q_j^{-1}(y_{j}')$ for $j=1,2$. Since $p$ is finite and
$x'\in X'$ is a general point we have $\rho_j
(p^{-1}(x'))=q_j^{-1} (x_{j}')$ (otherwise $\rho_j$ is not
dominant). Replacing $x$, if necessary, by another point from
$p^{-1} (x')$, we can suppose that $x_1=y_1$ and that $y_2=z_2$
for some $z \in p^{-1}(x')$ (this means that $y$ and $x$ belong to
the same orbit $O$ of $H_1$ because $x$ is general, see Remark
\ref{dimension}). Since $x'=z'$ we have $x_i'=y_i'=z_i'$ which
implies that $q_2$ sends $x_2$ and $z_2$ to the same point. By the
assumption $q_2$ is $H_1$-equivariant in case (i). In particular,
it sends the orbit $O_2=\rho_2 (O) \subset X_2$ into an
$H_1$-orbit $O_2' \subset X_2'$. Both orbits are isomorphic to
$H_1 \simeq \C_+$, i.e. the $H_1$-equivariant morphism
$q_2|_{O_2}: O_2 \to O_2'$ must be an isomorphism. That is,
$x_2=z_2=y_2$ and, therefore, $\rho (x)=(x_1,x_2)=(y_1,y_2)=\rho
(y)$. Since $\rho : X \to Z$ is birational and $x$ is general we
have $x=y$. Hence $\rho'$ is an embedding in a neighborhood of a
general point $x'$ which implies the desired conclusion for (i).

In case (ii) the general $\Gamma$-orbit $U$ in $X$ containing $O$
(resp. $\Gamma$-orbit $U' =p (U)\subset X'$ containing $O'=p(O)$)
is the set of left cosets of a finite subgroup $K$ (resp. $K'$) in
$SL_2 (\C )$. The $SL_2(\C )$-action is generated by
multiplication on the left while the $K$-action on $SL_2(\C )$ is
given by multiplication on the right. Hence the action of $\C_+
\simeq H_i < SL_2 (\C )$ commutes with the $K$-action. This
implies that each nonidentical element of $K$ sends any
$H_i$-orbit isomorphically into a different orbit.
Thus the quotient morphism $SL_2 (\C ) \to U$ (resp. $SL_2 (\C )
\to U'$) maps any $H_i$-orbit into a similar orbit isomorphically,
i. e. $p|_O : O \to O'$ is an isomorphism. If $U \simeq SL_2 (\C
)$ one can suppose that the restrictions of $\delta_1$ and
$\delta_2$ to $U$ are as in Example \ref{ex3} (2). The explicit
form of these derivations implies that $\rho_2 |_O : O \to O_2
=\rho_2 (O)$ is an isomorphism\footnote{Indeed, $x$ in this case
can be treated as a matrix $(c_{ij}) \in SL_2(\C )$, orbit $O$
consists of matrices $(c_{ij}(t))$ with $c_{1j}(t)=c_{1j}
+tc_{2j}$ and $c_{2j}(t)=c_{2j}$, while $\rho_2$ sends
$(c_{ij}(t))$ to vector $(c_{11}(t),c_{12}(t))$.}. Since the
$H_i$-actions on $SL_2(\C )$ commute with the $K$-action the same
is true in general case. Similarly, $\rho_2' |_{O'} : O' \to O_2'
=\rho_2' (O')$ is an isomorphism. Hence $q_2|_{O_2} : O_2 \to
O_2'$ is an isomorphism and the same argument as in case (i)
concludes case (ii).

\end{proof}

Though Lemma \ref{finereformulation} (together with a weak version of
Corollary \ref{holomorphic} that follows from it) enables us to go directly to
the proof of Theorem \ref{alg.groups}, we include some other
results to provide a stronger tool for establishing compatibility
condition.

\begin{lemma}\label{quasi} Let the assumption of Lemma
\ref{finereformulation} hold with one exception: instead of the
finiteness of $p$ we suppose that there are a surjective
quasi-finite morphism $r : S \to S'$ of normal affine algebraic
varieties equipped with trivial $\Gamma$-actions and a surjective
$\Gamma$-equivariant morphism $\varrho' : X' \to S'$ such that $X$
is isomorphic to fibred product $X' \times_{S'}S$ with $p : X\to
X'$ being the natural projection (i.e. $p$ is surjective
quasi-finite). Then the conclusion of Lemma
\ref{finereformulation} remains valid.

\end{lemma}

\begin{proof}

By the Noether normalization theorem, taking the spectrum of the
integral closure of $\C [S']$ in the field of rational functions
on $S$, we obtain a normal affine algebraic variety ${\tilde S}
\supset S$ with a finite morphism $\tr : \tS \to S'$ extending $r:
S\to S'$. Set $\tX = X' \times_{S'} \tS$ and denote by ${\tilde p}
: {\tilde X} \to X'$ the natural proejction. Then ${\tilde X} $
contains $X$ as a Zariski dense open subset, ${\tilde p}$ extends
$p$, and the $\Gamma$-action can be extended to ${\tilde X}$. Let
${\tilde \rho}_i : {\tilde X} \to {\tilde X}_i$ be the quotient
morphism of the $H_i$-action on $\tX$. For any nontrivial $f \in
\C [\tS ]$ whose zero locus contains $\tS \setminus S$ the
$f$-localizations (i.e. the localizations with respect to the
multiplicative system generated by $f$) of algebras $\C [ \tX ]$
and $\C [X]$ are isomorphic. Hence the similar localizations of
$\C [ \tX_i ]$ and $\C [X_i]$ are isomorphic. This implies that
the natural morphism $X_i \to \tX_i$ is an embedding (over $S
\hookrightarrow \tS$) and we have the following commutative
diagram
\[ \begin{array}{ccccccc}
X & \hookrightarrow & {\tilde X} & \stackrel{{\tilde
\rho}}{\rightarrow} & {\tilde X}_1\times {\tilde X}_2 &
\stackrel{{\tilde \tau}}{\rightarrow} & \tS\\
\, \, \, \, \downarrow {p} &&  \, \, \, \, \downarrow {\tilde p} &&
\, \,  \, \, \downarrow {\tilde q} && \, \,  \, \, \downarrow {\tilde r}\\
X' & = & X' & \stackrel{\rho'}{\rightarrow} &
X_1'\times X_2'& \stackrel{{\tau'}}{\rightarrow} & S'\\
\end{array} \]
where $\trho = (\trho_1 , \trho_2 )$ extends $\rho =(\rho_1,\rho_2) : X \to X_1
\times X_2$, $\varrho' = \tau' \circ \rho'$,
and the morphisms $\trho , \tq , {\tilde r}$ are finite.

Set $Z$ (resp. $Z'$, resp. $\tZ$) equal to the closure of $\rho
(X)$ in $X_1\times X_2$ (resp.  $\rho' (X')$ in $X_1'\times X_2'$,
resp.  $\trho (\tX )$ in $\tX_1\times \tX_2$).
By Proposition \ref{reformulation} $\rho : X \to Z$ is birational and hence
$\trho : \tX \to {\tilde Z}$ is birational being the extension of $\rho$. By Claim
in the proof of Lemma \ref{finereformulation} morphism $\rho' : X'
\to Z'$ is birational.
Note also that $\rho' : X' \to Z'$ is quasi-finite (indeed,
otherwise the commutative diagram implies that contrary to
Proposition \ref{reformulation} $\rho$ would not be quasi-finite
because $\tq$ is finite and $p$ is surjective).

Suppose that $z' \in Z'$ and $s'=\tau' (z')$. Since $r$ is
surjective one can choose $s \in S \subset \tS$ with $\tr (s)=s'$. Take $z \in
{\tilde q}^{-1} (z') \cap \rho (X)$ so that ${\tilde \tau} (z)=s$ (we can
do this because the natural projection $X \to S$ is surjective).
Hence $z' \in \rho' (X')$, i.e. $\rho'$ is surjective.
Furthermore, for any sequence $\{ x_i' \}$ of points in $X'$ such
that $\rho' (x_i') \to z'$ we can choose $x_i \in p^{-1}(x_i')$ so
that $\rho (x_i ) \to z$. Since morphism $\rho : X \to Z$ is
finite by Proposition \ref{reformulation}, one can suppose that
the sequence $\{ x_i \}$ is convergent to a point $x \in X$. Hence
the sequence $\{ x_i' \}$ is convergent to $x'=p (x)$ which means
that $\rho' : X' \to Z'$ is proper. Being also quasi-finite, this
morphism is finite by Grothendieck's theorem. Now Proposition
\ref{reformulation} yields the desired conclusion.

\end{proof}

\begin{remark}\label{gen.finite} We do not know whether the assumption, that
$\varrho'$ and $r$ are surjective, is essential. Without this
assumption the statement of Lemma \ref{quasi} says only that  $p$
is quasi-finite since one can put $S=X//\Gamma$ and
$S'=X'//\Gamma$. However the surjectivity of $p$ may be sufficient
for our purposes. Indeed, our aim is to check preservation of the
algebraic density property under quasi-finite morphisms and there
are examples of affine algebraic manifolds that are not
homogeneous (and, therefore, have no algebraic density property
\cite{V2}) but contain Zariski dense affine algebraic subvarieties
with the algebraic density property. For instance, the
hypersurface in $\C^3_{x,y,z}$ given by $xy=z^2-1$ has the
algebraic density (\cite{KK}) and it is not difficult to show that
it is isomorphic to the complement to the line $x=z=1$ in the
hypersurface in $\C^3$ given by $x(x-1)y=z^2-1$. The dual graph of
a simple normal crossing completion of the latter hypersurface
cannot be contracted to a zigzag in the terminology of \cite{Gi}
and, hence, this hypersurface is not even quasi-homogeneous by
Gizatullin's theorem.

\end{remark}

Recall that an \'etale neighborhood of a point $y$ of an algebraic
variety $Y$ is an \'etale morphism $g : W \to Y$
whose image contains $y$.

\begin{proposition}\label{condition2} Let $Y$ be a normal
affine algebraic variety equipped with a trivial $\Gamma$-action
(where $\Gamma$ is from Lemma \ref{finereformulation}) and $r : X
\to Y$ be a surjective $\Gamma$-equivariant morphism. Suppose that
for any $y \in Y$ there exists an \'etale neighborhood $g: W \to
Y$ such that the vector fields induced by $\delta_1$ and
$\delta_2$ on the fibred product $X\times_Y W$ are compatible.
Then $\delta_1$ and $\delta_2$ are compatible.
\end{proposition}

\begin{proof} Set $Y_1=g(W)$. Then the restrictions $\delta_1^1$
and $\delta_2^1$ of $\delta_1$ and $\delta_2$ to $X_1=r^{-1}(Y_1)$
are compatible by Lemma \ref{quasi}. Suppose that $\{ Y_i \}$ is a
finite cover of $Y$ by open sets similar to $Y_1$ and notation
$X_i, \delta_1^i , \delta_2^i$ have also the similar meaning.
Without loss of generality we can assume that $Y_i =Y \setminus
f_i^{-1}(0)$ for some $f_i \in \C [Y_i] \subset \C [X_i]$. Let
$I_i\subset \C [X_1]$ be the largest ideal contained in ${\rm
Span}(\Ker \delta_1^i \cdot \Ker \delta_2^i)$ and $I$ be the
largest ideal in $\C [X]$ whose $f_i$-localization is contained in
$I_i$ for every $i$. In particular, $I$ is non-zero since each
$I_i$ is such. Show that $I \subset {\rm Span}(\Ker \delta_1 \cdot
\Ker \delta_2)$.

Indeed, $f_i \in \C [Y] \subset \Ker \delta_j, \, j=1,2$. Hence
for every $a \in I$ there exists $k_i$ such that $af_i^{k_i}$ is in ${\rm Span}(\Ker \delta_1
\cdot \Ker \delta_2)$. By Hilbert's Nullstellensatz there are
regular functions $g_i$ on $Y$ such that $\sum_if_i^{k_i}g_i
\equiv 1$. Since $g_i$ is in the kernel of $\delta_1$ we see that
$a \in {\rm Span}(\Ker \delta_1 \cdot \Ker \delta_2)$ which
concludes the proof.

\end{proof}

\begin{corollary}\label{holomorphic}
Let a linear algebraic group $G$ act algebraically  on $X$ so that
$X//G$ is affine, the quotient morphism $X \to X//G$ is surjective
(which is always true when $G$ is reductive) and makes $X$ an
\'etale $G$-principal bundle over $X//G$. Suppose that $\Gamma$
(from Lemma \ref{finereformulation}) is an algebraic subgroup of
$G$ and the actions of $H_i, \, i=1,2$ on $G$ induced by left
multiplication generate compatible derivations on $\C [G]$. Let
the induced $H_i$-actions on $X$ correspond to derivations
$\delta_i$ on $\C [X]$. Then $\delta_1$ and $\delta_2$ are
compatible.
\end{corollary}

\begin{theorem}\label{alg.groups} Let $G$ be a linear algebraic group whose connected
component is different from a torus or $\C_+$. Then $G$ has the
algebraic density property.

\end{theorem}

\begin{proof}
Since all components of $G$ are isomorphic as varieties we can
suppose that $G$ is connected. Recall that the unipotent radical
$R$ of $G$ is an algebraic subgroup of $G$ (\cite{Ch}, p. 183). By
Mostow's theorem \cite{Mo} (see also \cite{Ch}, p. 181) $G$
contains a (Levi) maximal closed reductive algebraic subgroup $L$
(which is, in particular, affine) such that $G$ is the semi-direct
product of $L$ and $R$, i.e. $G$ is isomorphic as affine variety
to the product $R \times L$. In case $L$ is trivial $G=R \simeq
\C^n, n\ge 2$ and  we are done by Corollary \ref{AL}. In the case
of both $R$ and $L$ being nontrivial we are done by Corollary
\ref{product} with $R$ playing the role of $X_1$ and $L$ of $X_2$.

Thus it remains to cope with reductive groups $G$. Let $Z\simeq
(\C^*)^n$ denote the center of $G$ and $S$ its semisimple part.
First we suppose that $Z$ is nontrivial. The  case when $G$ is
isomorphic as group to the direct product $S\times Z$  can be
handled as above by Corollary \ref{product} with $S$ playing the
role of $X_1$ and $Z$ of $X_2$. In particular, we have a finite
set of pairs of compatible vector fields $\{ \delta_1^k ,
\delta_2^k \}$ as in Theorem \ref{density}. Furthermore, one can
suppose that the fields $\delta_1^k$ correspond to one parameter
subgroups of $S$ isomorphic to $\C_+$ and $\delta_2^k$ to one
parameter subgroups of $Z$ isomorphic to $\C^*$. In the general
case $G$ is the factor group of $S\times Z$ by a finite (central)
normal subgroup $\Lambda$. Since $\Lambda$ is central the fields $
\delta_1^k$, $\delta_2^k $ induce completely integrable vector
fields $\tilde \delta_1^k$, $\tilde \delta_2^k $ on $G$ while
$\tilde \delta_2^k (x_0)$ is a generating set for some $x_0 \in
G$. By Lemma \ref{finereformulation} the pairs $\{ \tilde
\delta_1^k, \tilde \delta_2^k \}$ are compatible and the density
property for $G$ follows again from Theorem \ref{density}.


It remains now to consider a semi-simple $G$ which can be assumed
simply connected by Lemma \ref{finereformulation}. That is, it is
a product of simple Lie groups and by Corollary \ref{product} it
suffices to consider the case when $G$ is simple. Such $G$
contains $SL_2(\C )$ as a subgroup.
%
%
%
%
The existence of two compatible vector fields $\delta_1$ and
$\delta_2$ on $SL_2( \C )$ implies their existence on $G$ by
Corollary \ref{holomorphic}. Since a simple Lie group is
tangentially semi-homogenous (see Example \ref{ex3}) the algebraic
density property for $G$ follows again from Theorem \ref{density}.

\end{proof}

\section{Codimension 2 Case.}

{\bf Motivation and Notation.} In this section $X$ will be a
closed affine algebraic subvariety of $\C^n$ whose codimension
$n-k$ is at least 2. By the Hartogs theorem any completely
integrable algebraic (or holomorphic) vector field  on $\C^n
\setminus X$ extends to a similar vector field on $\C^n$ tangent
to $X$. In particular, the Lie algebra generated by completely
integrable algebraic (or holomorphic) vector fields contains only
vector fields tangent to $X$, i.e. there is no density property
for $\C^n \setminus X$. In general there is no also hope that this
Lie algebra coincides with the Lie algebra of all algebraic vector
fields tangent to X, since this would imply density property for
$X$ (and our $X$ maybe even not smooth!). Therefore, it is natural
to study the Lie algebra $\LieA (\C^n, X)$ generated by completely
integrable algebraic vector fields on $\C^n$ that vanishes on $X$.
According to Forstneri\v c the best possible result to expect is
that $\LieA (\C^n, X)$ is equal to the Lie algebra of all
algebraic vector fields vanishing on $X$. We use notation $\AVF_I
(\C^n)$ for the latter algebra where $I\subset \C^{[n]}$ is the
defining ideal of $X$ (more generally, for any affine algebraic
variety $Y$ and an ideal $L \subset \C [Y]$ we denote by
$\AVF_L(Y)$ the Lie algebra of vector fields whose coordinate
functions are from $L$). If the above property holds then the
geometric structure of algebraic vector fields vanishing on $X$
has the algebraic density property in the terminology of Varolin
\cite{V2}.

We will prove this property under some weak additional assumption
in Theorem \ref{codimension} and a very close result without any
additional assumption in Theorem \ref{multiple}. Both results lead
to a generalization of the main theorem of the Anders\'en- Lempert
theory, allowing now to construct holomorphic automorphisms of
$\C^n$ not only with control on compacts but with additional
control on algebraic subvarieties of codimension at least 2.

\begin{lemma}\label{GW}   The group
$\Aut (\C^n, X)$ of algebraic  automorphisms of $\C^n$ identical on $X$
acts transitively on $Z=\C^n\setminus X$ and, furthermore, for any $z \in Z$ the image  of any vector $v\in T_z Z$ under
the isotropy group $\Aut (\C^n, X) _z$ generates $T_z Z$ (compare with definition \ref{tangential}).
\end{lemma}

\begin{proof}
By a theorem of Gromov \cite{Gro} and Winkelmann \cite{Wi} $Z$ is
homogenous. We will use the idea of their proof.
More precisely, consider a general linear projection $p: \C^n \to
\cH \simeq \C^{n-1}$ and a nonzero constant vector field $\nu$
such that $p_* (\nu )=0$. Then $p(X)$ is a subvariety of codimension
at least 1 in $\cH$. For every regular function $h$ on $\cH$ that
vanishes on $p(X)$ the vector field $h \nu$ generates a
$\C_+$-action on $Z$. Changing $\cH$ we get a transitive action.

Consider a general point $z \in Z$ whose projection $z_0 \in \cH$
is not in $p(X)$. Suppose that $h$ has a simple zero at $z_0$. By
the Claim in the proof of Corollary \ref{product} the
$\C_+$-action generated by $h \nu$ acts on $T_zZ$ by  the formula
$w \to w+ dh (w) \nu (w)$ where $dh$ is the differential of $h$
and $w \in T_zZ$. Since $\nu$ may be chosen as a general constant
vector field on $\C^n$ we see that $\Aut (\C^n , X)_z$ induces an
irreducible representation on $T_zZ$ which implies the second
statement.

\end{proof}

\begin{theorem}\label{multiple}
There is an ideal $L\subset \C^{[n]}$ whose radical is $I$ such
that $\LieA (\C^n, X)$ contains $\AVF_L (\C^n)$.

\end{theorem}

\begin{proof} Suppose that $x_1, \ldots ,x_n$ is a coordinate
system, $p_i: \C^n \to \C^{n-1}$ is a projection to the coordinate
hyperplane $\cH_i =\{ x_i = 0 \}$, and $h_i$ is a nonzero function
on $\cH_i$ that vanishes on $p_i(X)$. Set $\delta_i =
\partial /\partial x_i$ and choose $f_i \in \Ker \delta_i$. Then
$f_ih_i\delta_i$ is a completely integrable algebraic vector field
on $\C^n$ that vanishes on $X$, i.e. it generates a $\C_+$-action
on $Z$ (since the elements of this action are from $\Aut (\C^n, X)$). Then
$$[f_1h_1\delta_1,x_1f_2h_2\delta_2]-[x_1f_1h_1\delta_1,f_2h_2\delta_2]=f_1f_2h_1h_2\delta_2$$
belongs to $\LieA (\C^n, X)$. Since $\Ker \delta_1 \cdot \Ker \delta_2$
generates the ring of polynomials $\C^{[n]}$ as a vector space we
see that $\LieA (\C^n, X)$ contains all algebraic fields proportional to
$\delta_2$ with coordinate functions in the principal ideal
generated by $h_1h_2$. Since one can perturb $x_2$ (as a linear function)
 $\LieA (\C^n, X)$ contains all algebraic vector fields whose coordinates are in some (non-zero)
 ideal $L$. Since $Z$ is homogenous under $\Aut (\C^n, X)$ arguing as in the proof of
 Theorem \ref{semi} one can suppose that the radical of $L$ is $I$.

\end{proof}

Though Theorem \ref{multiple} does not give the algebraic density
of the Lie algebra of algebraic vector fields vanishing on $X$ it
implies already a strong approximation result generalizing the
Anders\' en-Lempert theorem. We omit its proof since it repeats
the arguments in \cite{FR} with minor modifications.

\begin{theorem}\label{approximation} Let $X$ be an algebraic subvariety of $\C^n$ of
codimension at least 2
and $\Omega$ be an open set in $\C^n$ $(n\ge
2)$. Let $\Phi : [0, 1]\times \Omega \to \C^n \setminus X$ be a
${\cC}^2$-map
such that for every $t \in [0, 1]$ the restriction $\Phi_t=\Phi |_{(t
\times \Omega )}$ is an injective holomorphic map onto a Runge domain
$\Phi_t (\Omega)$ in
$\C^n$. If $\Phi_0$ can be approximated on $\Omega$ by holomorphic
automorphisms of
$\C^n$  identical on $X$, then for every $t \in [0, 1]$ the map $\Phi_t$
can be approximated on $\Omega$ by such automorphisms.
\label{ALTH}
\end{theorem}

As a consequence we recover the result of Buzzard
and Hubbard \cite{BH} answering Siu's question.

\begin{corollary} Any point $z$ in the complement of an algebraic subset $X$ of $\C^n$ of codimension at least 2
has a neighborhood $U$ in $\C^n\setminus X$ that is biholomorphic
to $\C^n$ (such U is called a Fatou-Bieberbach domain).
\end{corollary}

\begin{proof} Following the standard scheme of Rudin and Rosay,
choose a ball $\Omega \subset \C^n \setminus X$  around $z$ and consider $\Phi_t$
contracting this ball radially
towards $z$. The resulting automorphism approximating $\Phi_1$ from Theorem \ref{ALTH}  has an attracting
fixed point near $z$ and $z$ is contained in the basin of attraction. This basin is a Fatou-Bieberbach domain
and it does not meet $X$ since the automorphism is identical on $X$.
\end{proof}

Though Theorem \ref{multiple} is sufficient for the approximation
purposes, let us be accurate and establish the algebraic density
for algebraic vector fields vanishing on $X$ under an additional
assumption.

\begin{convention} We suppose further in this section that the dimension of
the Zariski tangent space $T_xX$ is at most $n-1$ for every point
$x \in X$.

\end{convention}

\begin{lemma}\label{technical} Lie algebra $\LieA (\C^n, X)$ contains $\AVF_{I^2} (\C^n )$.
\end{lemma}

\begin{proof} It suffices to show that for every point $o \in \C^n$ there
exists a Zariski neighborhood $V$ and a submodule $M_V$ from
$\LieA (\C^n, X)$ such that its localization to $V$ coincides with
the localization of $\AVF_{I^2} (\C^n )$ to $V$. Indeed, because
of quasi-compactness we can find a finite number of such open sets
$V_i$ that covers $\C^n$. Hence the coherent sheaves generated by
$\AVF_{I^2} (\C^n )$ and $\sum_i M_{V_i}$ coincide locally which
implies that they have the same global sections over affine
varieties by Serre's theorem B. In fact,  it suffices to show that
the localization of $M_V$ to $V$ contains all fields from the
localization of $\AVF_{I^2} (\C^n )$ to $V$ that are proportional
to some general constant vector field $\delta$ which is {\em our
aim now}. By Theorem \ref{multiple} it is also enough to consider
$o \in X$ only. The construction of the desired neighborhood $V$
of $o$ starts with the following.

{\em Claim.} For any point $o \in X$,
$l \geq \max (k+1,\dim
T_oX)$ (where $\dim X =k$), and a general linear projection $p: \C^n \to \cH \simeq
\C^l$ one can choose a projection $p_0 : \C^n \to \cH_0 \simeq
\C^{l-1}$ for which

{\rm (i)} $p_0=\varrho \circ p$ where $\varrho : \cH \to \cH_0$ is
a general linear projection, and

{\rm (ii)} there exists $h \in \C [\cH_0] \simeq \C^{[l-1]} \simeq
\varrho^*(\C^{[l-1]}) \subset \C^{[l]}$ such that $h$ does not
vanish at $p_0(o)$ and $p|_{X \setminus (h \circ p)^{-1}(0)} : X
\setminus (h \circ p)^{-1}(0) \to p(X) \setminus h^{-1}(0)$ is an
isomorphism.

Since $p$ is general the condition on $l$ implies that $p$ is a
local isomorphism in a neighborhood of $o$ and, furthermore, since
$\varrho$ is also general then by Bertini's theorem
$p_0^{-1}(p_0(o))$ contains only smooth points of $X$ except, may
be, for $o$, i.e. $p$ is a local isomorphism in a neighborhood of
each of these points which implies the Claim.

From now on let $l=n-1$.
Choose a general coordinate system  ${\bar x}=(x_1, \ldots ,x_n)$
on $\C^n$ such that $p({\bar x})=(x_2, \ldots ,x_n)$ and
$p_0({\bar x})=(x_3, \ldots ,x_n)$, i.e. $h=h(x_3, \ldots ,x_n)$.
Set $V=\C^n \setminus h^{-1}(0)$.

Since $p(X) \cap V \simeq X \cap V$ we have $x_1=r/h^s$ where $r$
is a polynomial in $x_2, \ldots , x_n$ and $s \geq 0$. Set
$\nu_i=\partial /\partial x_i$ for $i \ne 2$, and $\nu_2=h^s
\partial /\partial x_2 +(\partial r/\partial x_2) \partial
/\partial x_1$. Then each $\nu_i$ is a locally nilpotent
derivation and $\Ker \nu_1$ contains the defining ideal $I_p$ of
$p(X)$ in $\C [\cH ] \simeq \C^{[n-1]} \simeq p^*(\C^{[n-1]})
\subset \C^{[n]}$. Furthermore, for $\xi =h^sx_1 -r$ we have $\xi
\in \Ker \nu_2$, and $\xi$ (resp. $x_2$) is of degree 1 with
respect to $\nu_1$ (resp. $\nu_2$). This implies that for $f,g \in
I_p$ the vector fields that appear in the Lie brackets below are
completely integrable and vanish on $X$:

$$ [f\nu_1,\xi g\nu_1]=h^sfg\nu_1 , \, \,
[\xi \nu_2 , x_2 \xi \nu_1] -[x_2 \xi \nu_2 , \xi \nu_1]=h^s\xi^2
\nu_1 ,$$   $${\rm and} \, \, \, [\xi \nu_2 , x_2f \nu_1]-[x_2 \xi
\nu_2 , f \nu_1 ]=h^s\xi f \nu_1. $$

The defining ideal of $X \cap V$ is generated by $\xi$ and
elements of $I_p$. Since $h$ is invertible on $V$ and $\nu_1$ is a general
constant vector field
from the formulas
before we see that the localization of $\LieA (Z)$ to $V$ contains
the localization of $\AVF_{I^2} (\C^n ) $ which is
the desired conclusion.

\end{proof}

\begin{theorem}\label{codimension}
Let $X$ be a closed algebraic subset of $\C^n$ of codimension at
least 2 such that the Zariski tangent space $T_xX$ has dimension
at most $n-1$ for any point $x \in X$. Then $\LieA (\C^n, X)$ is equal
to  $\AVF_I (\C^n)$, i.e. the vector fields  vanishing on $X$ have algebraic density property.\end{theorem}

\begin{proof} Similarly to the proof of Lemma \ref{technical}, it suffices
to show that for every point $o \in \C^n$ there exists a Zariski
neighborhood $V$ and a submodule from $\LieA (\C^n, X )$ such that
its localization $M$ to $V$ coincides with localization of $\AVF_I
(\C^n)$ to $V$. By Theorem \ref{multiple} it is enough to consider
$o \in X$ and, furthermore, it suffices to show that this
localization $M$ contains all elements of $\AVF_I (\C^n)$
proportional to some general constant vector field.

Let $\nu_i, p, I_p$, and $\xi$ have the same meaning as in the
proof of Lemma \ref{technical}. Choose $\nu_1$ as this constant
vector field. Since $I$ is generated by $\xi$ and $I_p$ one needs
to show that all fields of the form $\mu =(\xi g_0 +\sum g_i f_i)
\nu_1$ are contained in $M$ where $g_0,g_i$ are regular on $V$ and
$f_i \in I_p$. Since $p$ yields an isomorphism between $p(X) \cap
V $ and $ X \cap V$ there are functions $e_0,e_i$ that do not
depend on $x_1$ and such that $e_0|_X =g_0|_X$ and
$e_i|_X=g_i|_X$. Then $\mu =(\xi e_0 +\sum e_i f_i) \nu_1 + a
\nu_1$ where $a$ belongs to the localization of $I^2$ to $V$ (e.g.
$a= \xi( g_0-e_0) +\sum (g_i -e_i) f_i$). Since the first summand
in the last formula for $\mu$ is completely integrable we have the
desired conclusion from Lemma \ref{technical}.

\end{proof}

\begin{remark}

(1) The authors believe that the condition $\dim T_xX \le n-1$ in
Theorem \ref{codimension} is essential. As a potential
counterexample one may try to take $X$ equal to polynomial curve
in $\C^3$ with one singular point whose Zariski tangent space is
$3$-dimensional. More precisely, let $L$ be the Lie algebra
generated by vector fields that vanish on $X$ and have form $f
\sigma$ where $f \in \Ker \sigma$ and $\sigma$ is either locally
nilpotent or semi-simple. Then we can show that $L$ does not
coincide with $\AVF_I (\C^3)$ but we do not know whether $L=\LieA
(\C^3,X)$.

(2) In view of Theorem \ref{codimension} the assumptions of Theorem \ref{ALTH} can
be weakened in case of $\dim T_xX \le n-1$  to the following extend:
the assumption $\Phi_t (\Omega)\cap X = \emptyset$  can be replaced by
the assumption that $\Phi_t$ is identical on
$\Phi_t (\Omega)\cap X$ for all $t$.

(3) The assumption of codimension at least 2 for $X$ cannot be
removed, since the complement to a hypersurface in $\C^n$ can be of general type or even
Kobayashi hyperbolic and, therefore, there is no nontrivial
completely integrable holomorphic vector field vanishing on $X$.
Also the assumption that $X$ is not just a holomorphic but an
algebraic subvariety of $\C^n$ cannot be weakened. This follows
from the fact that there are holomorphic embeddings of $\C$ into
$\C^n$ (for any $n$) such that the group of holomorphic
automorphisms of $\C^n$ identical on the image is trivial (e.g.,
see \cite{DKW}).
\end{remark}

\providecommand{\bysame}{\leavevmode\hboxto3em{\hrulefill}\thinspace}

\end{document}